\documentclass[12pt,a4paper,reqno]{amsart} 
\usepackage{tikz}
\usepackage[utf8]{inputenc}
\usepackage{amsmath,bm,bbm,amsthm, amssymb}
\usepackage{color}
\usepackage{animate}
\usepackage{etoolbox}
\usepackage[ruled,vlined,linesnumbered,algosection,resetcount]{algorithm2e}
\usepackage{comment}
\usepackage{pgf}
\usepackage{hyperref}
\usetikzlibrary{calc}
\usetikzlibrary{patterns}
\usetikzlibrary{arrows}
\usetikzlibrary{arrows.meta}
\usetikzlibrary{decorations.pathreplacing}
\usepackage[utf8]{inputenc}
\usepackage{dsfont}

\theoremstyle{plain}
\newtheorem{theorem}{Theorem}
\newtheorem{proposition}[theorem]{Proposition}
\newtheorem{corollary}[theorem]{Corollary}
\newtheorem{lemma}[theorem]{Lemma}

\theoremstyle{definition}
\newtheorem{definition}[theorem]{Definition}

\newtheorem{example}[theorem]{Example}

\theoremstyle{remark}

\definecolor{wiasblue}   {cmyk}{1.0, 0.60, 0, 0}
\definecolor{mlugreen}{RGB}{0,81,51}

\def\Z{\mathbb Z}
\def\E{\mathbb E}
\def\P{\mathbb P}

\def\R{\mathbb R}
\def\mc{\mathcal}
\def\ms{\mathsf}

\newcommand{\pimsta}[1]{\pi_{-s, t, T}}

\def\vp{\varphi}

\def\to{\uparrow}

\def\th{\theta}

\def\su{\subseteq}
\def\De{\Delta}

\def\Cov{\ms{Cov}}

\def\been{\begin{enumerate}}
\def\enen{\end{enumerate}}
\def\im{\item}

	\def\bec{\begin{corollary}}
	\def\enc{\end{corollary}}

\def\tff{\to \ff}

\def\CC{\mathcal{C}}

\def\lrsa{\leftrightsquigarrow}
\def\la{\lambda}
\def\La{\Lambda}

\def\a{\alpha}
\def\s{\sigma}
\def\su{\subseteq}
\def\e{\varepsilon}
\def\t{\tau}
\def\g{\gamma}

\def\de{\delta}

\def\es{\emptyset}
\def\one{\mathbbmss{1}}

\def\De{\Delta}

\def\co{\colon}
\def\ff{\infty}

\def\lc{\lceil}
\def\rc{\rceil}

\def\vp{\varphi}

\def\d{{\rm d}}
\def\k{\kappa}

\def\FF{\mathcal{F}}

\def\lc{\lceil}
\def\rc{\rceil}
\def\lf{\lfloor}
\def\rf{\rfloor}

\def\f{\frac}

\def\xra{\xrightarrow}
\def\xrd{\xrightarrow{\hspace{.2cm}\ms D\hspace{.2cm}}}

\def\off{[0,\ff)}

\def\ua{\uparrow}

\def\en{\end}
\def\im{\item}
\def\sm{\setminus}
\def\th{\theta}
\def\bep{\begin{proof}}
\def\enp{\end{proof}}
\def\bepr{\begin{proposition}}
\def\enpr{\end{proposition}}
\def\bec{\begin{corollary}}
\def\enc{\end{corollary}}
\def\bea{\begin{align}}
\newcommand\eea{\end{align}}
\def\beas{\begin{align*}}
\def\eeas{\end{align*}}
\def\bet{\begin{theorem}}
\def\ent{\end{theorem}}
\def\bee{\begin{example}}
\def\ene{\end{example}}

\def\da{\downarrow}
\def\bede{\begin{definition}}
\def\ende{\end{definition}}
\def\bel{\begin{lemma}}
\def\enl{\end{lemma}}
\def\been{\begin{enumerate}}
\def\enen{\end{enumerate}}

\def\beit{\begin{itemize}}
\def\enit{\end{itemize}}
\def\befr{\begin{frame}}
\def\enfr{\end{frame}}

\def\ti{\times}

\def\Cov{\ms{Cov}}
\def\ba{\,|\,}

\def\pa{\partial}

\def\tf{\tfrac}
\def\ns{n_{\ms S}}
\def\nsp{n_{\ms S}'}
\def\bef{\begin{figure}[!h]}
\def\enf{\end{figure}}
\def\bs{\boldsymbol}

\def\ms{\mathsf}
\def\mb{\mathbb}

\def\Cov{\ms{Cov}}

\def\ccn{\mc C_{\de_0, \ff}}

\def\tktb{\t_T}

\def\vli{\lambda_{\ms S}}

\def\d{{\rm d}}
\def\g{\gamma}
\def\P{\mb{P}}
\def\es{\emptyset}
\def\s{\sigma}
\def\la{\lambda}
\def\a{\alpha}

\def\eq{\begin{equation}}
\def\en{\end{equation}}

\def\itf{[0, 1] \rightarrow [0, 1]}
\def\jtf{[0, \ff) \rightarrow [0, 1]}

\def\idm{I_{\de, M}}
\def\idmm{I_{k, L, \de, M}}
\def\ildm{I_{o, L, \de, M}}
\def\jdm{J_{\de, M}}

\usepackage{setspace}
\calclayout

\keywords{wireless communication network,  multi-scale model, scaling limit, continuum percolation}
\subjclass[2010]{60K35; 60F10; 82C22}
\date{\today}

\begin{document}

\title{Connection intervals in multi-scale dynamic networks}

\author{Christian Hirsch}
\address[Christian Hirsch]{Bernoulli Institute, University of Groningen, Nijenborgh 9, 9747 AG Groningen, The Netherlands}
\email{c.p.hirsch@rug.nl}
\author{Benedikt Jahnel}
\address[Benedikt Jahnel]{Weierstrass Institute for Applied Analysis and Stochastics, Mohrenstra\ss e 39, 10117 Berlin, Germany}
\email{benedikt.jahnel@wias-berlin.de}
\author{Elie Cali}
\address[Elie Cali]{Orange SA, 44 Avenue de la R\'epublique, 92326 Ch\^atillon, France}
\email{elie.cali@orange.com}

\begin{abstract}
We consider a hybrid spatial communication system in which mobile nodes can connect to static sinks in a bounded number of intermediate relaying hops. We describe the distribution of the connection intervals of a typical mobile node, i.e., the intervals of uninterrupted connection to the family of sinks. This is achieved  in the limit of many hops, sparse sinks and growing time horizons. We identify three regimes illustrating that the limiting distribution depends sensitively on the scaling of the time horizon. 
\end{abstract}

\maketitle

\section{Introduction}
\label{int_sec}
Starting with the landmark paper~\cite{gilbert} in the early 1960s by Gilbert, {stochastic geometry} has been employed to model and analyze spatial communication systems in which the network nodes directly exchange data with other nodes in their vicinity. In the absence of any refined information about the spatial locations of nodes the null model is that they are scattered entirely at random in space, i.e., form a \emph{homogeneous Poisson point process}, where the single scalar parameter represents the expected number of vertices per unit volume. Concerning the communication structure, the simplest model is the \emph{Gilbert graph} where connections between nodes are represented by links between any pair of nodes with a certain maximal distance. 

Basic questions about the connectivity of such {peer-to-peer} networks provide a key motivation for fruitful research in the realm of continuum percolation. At the center of this field stands the {percolation phase transition}, meaning that if the intensity of network participants is sufficiently high, then a positive proportion of all nodes form a giant communicating cluster \cite{cPerc}.


However, only a small set of use cases such as {sensor networks} or {disaster-rescue ad-hoc network} rely on peer-to-peer networks in its purest sense. In the bulk of applications, peer-to-peer communications appears as an extension for more traditional {cellular networks}, forming a variety of {hybrid systems}~\cite{mscale}. Such systems have the potential to successfully mitigate many of the problems of the pure systems, such as for example {delay}, {jitter}, {routing} or {operational control}.

Another essential aspect is {mobility}. The vast majority of the available literature on stochastic models for spatial communication networks investigate static systems. However, already in the landmark paper~\cite{grossTse}, the impact of mobility on the capacity of communication networks has been evaluated in an information-theoretic context. These findings have inspired subsequent studies of spatial random networks with mobility, and we refer the reader to \cite{diMi,adHoc} for an overview in this area.

When designing and evaluating hybrid communication networks, arguably the most basic network characteristic is the \emph{total connection time}. That is, the overall time that a typical network node is connected to some infrastructure. The key achievement of our earlier work \cite{mscale} is to describe the asymptotic behavior of this quantity over long time horizons, many hops and sparse infrastructure nodes. However, in real communication networks, maximizing the total connection time does not necessarily lead to networks offering an acceptable quality of service. Indeed, if the connection times are highly fragmented over the entire time horizon, then it is not possible to offer the typical node a large coherent block of uninterrupted service, and the system faces a substantial overhead caused by the cost of frequently re-establishing lost connections to the typical node. Therefore, in the present work, we move beyond the total connection time and provide more refined descriptions of the \emph{connection intervals}, i.e., the time intervals when the typical mobile node is guaranteed uninterrupted connection. In particular, this distribution can be used to answer questions of the following form.

\begin{enumerate}
\item What is the proportion of time that a typical mobile node is guaranteed uninterrupted communication of at least a given time duration?
\item What is the number of reconnections of a typical mobile node in the hybrid system?
\end{enumerate}


The rest of the manuscript is organized as follows.  In Section~\ref{Sec_Set}, we introduce precisely the hybrid model and state Theorem~\ref{thm_2} as the main result on the weak convergence of the connection-interval measure under three different coupled limits. We also present a simulation study to illustrate how the results can be applied for designing wireless networks. In Section~\ref{out_sec}, we outline the proof of Theorem \ref{thm_2} and establish a number of supporting results that feature several approximations. Finally, Section~\ref{Sec_Pro} contains the detailed proofs of the supporting results.

\section{Setting and main result}
\label{Sec_Set}

\subsection{System model}
\label{sys_sec}

In this manuscript, we study an infrastructure-enhanced model for a wireless communication network in $\R^d$, which is observed over a time horizon $[0, T]$. More precisely, the infrastructure nodes or sinks form a homogeneous Poisson point process $Y = \{Y_j\}_{j \ge 1}$ with some intensity $\vli > 0$.

Furthermore, at every time instant $t \le T$ also the mobile nodes $X(t)=\{X_i(t)\}_{i \ge 1}$ form a homogeneous Poisson point process with intensity $\la > 0$. The nodes move over time according to a mobility model specified further below.

We work with a simple distance-based connection model, where nodes and sinks communicate directly if they are closer than a certain fixed communication radius. By the scaling properties of the Poisson point process, we may henceforth fix this distance to be equal to 1.  Moreover, also the nodes can communicate among themselves at this distance and they may forward data over several relaying hops. Figure~\ref{net_fig} (left) illustrates such a network without any sinks. Thus, a $k$-hop connection with a sink can be established even if the latter is outside the range of direct connectivity.

\subsection{Connection intervals}\label{Sec_Co}

In this work, we go beyond the setting of empirical $k$-hop connection times as studied in~\cite{mscale} and investigate the connection intervals. More precisely, we let $I(t, S) \su \R$ denote the length of the connected component of a set $S \su \R$ containing a time point $t \in \R$. In symbols, 
$$
I(t, S):= \sup_{a \le b \co t \in [a, b] \su S}(b - a),
$$
where $I(t, S) := 0$ if $t \not \in S$.

In the specific setting of the network model from Section \ref{sys_sec}, $S$ represents the set of all times where a typical moving node $X_0(\cdot)$ started at the origin is $k$-hop connected to some sink.  Note that we may add this extra point, due to the Slivnyak theorem, and after applying mobility to this point as well as to the other points, it remains typical at every instant. More precisely, 
$$\Xi^k(Y_j) := \{t \in \R\co X_0(t) \stackrel k{\lrsa_t}  Y_j\}$$
denotes the set of all times when the typical node connects  to the sink $Y_j$ in at most $k$ hops using the nodes in $X$. Then, 
$$\Xi^k := \bigcup_{j \ge 1} \Xi^k(Y_j)$$
is the set where it connects to some sink. In particular, this definition allows for hand-overs between different sinks.

Writing $\de$ for the Dirac measure, the central objective of this work is the \emph{$k$-hop connection-interval measure} 
$$
\tktb(\d \ell, \d t) := \frac1T\int_{[0, T] \cap \Xi^k}\de_{(I(s, \Xi^k ), s/T)}(\d \ell, \d t) \d s,$$
which is a random probability measure on $[0, \ff) \ti[0, T]$. In order to highlight the richness of information encoded in the connection interval measure, we now describe three key network statistics derived from it.
\bee[Network statistics]\phantom{}
\label{stat_ex}
Consider $\t_T(f) := \int_{[0,\ff) \ti [0, T]} f(\ell, t) \t_T(\d \ell, \d t)$ for 
\been
\im {$\bs{ f_1(\ell, t) := 1.}$} Then, this is the time-averaged connection time of the typical node. Hence, we recover the setting from \cite{mscale}.
\im {$\bs{ f_2(\ell, t) := \ell.}$} As elucidated in Section \ref{int_sec}, connection times cannot be used effectively if they are highly fragmented over the time horizon. One approach is to use $f(\ell, t) = \one\{\ell \ge t_{\ms{min}}\}$, i.e., to discard connection times contained in intervals shorter than a minimum duration $t_{\ms{min}}$. However, this hard threshold might cause undesirable threshold phenomena tied to the specific choice of $t_{\ms{min}}$. Hence, it may be more desirable to rely on a soft weighting of the form $f_2(\ell, t) = \ell$, where connection times in longer intervals receive a higher weight. Then, the network characteristic $\t_T(f)$ takes into account that a network should not only offer a high amount of connectivity on average but also guarantee connections that are uninterrupted for a substantial time.
\im {$\bs{ f_3(\ell, t) := 1\{ \ell  > 0\}/\ell.}$} Then, $T\t_T(f)$ is the number of connection intervals in $[0, T]$ with the first and last such interval possibly counted only partially. Hence, $T\t_T(f)$ can be interpreted as the number of reconnections that are needed within the time horizon. For network operators this characteristic is of high interest since each such reconnection requires additional resources.
\enen
\ene

%
%
\subsection{Asymptotic $k$-hop connection-interval measure}
We analyze the connection measure $\tktb$ asymptotically over a growing time horizon $T \to \ff$ when simultaneously the admissible number of hops $k$ and the sink intensity $\vli$ scale with $T$. More precisely, 
$$\la_c := \inf\big\{\la > 0 \co \P(o \lrsa \ff) > 0\big\}$$ 
denotes the critical intensity for percolation of the node process. For $\la > \la_c$ we write $\CC_\ff = \CC_\ff^\la$ for the unique connected component in a Poisson point process with intensity $\la$.
We will rely on a key finding from continuum first-passage percolation, namely that above $\la_c$, the number of hops needed to travel inside the infinite connected component of nodes grows linearly in the Euclidean distance, see \cite{yao}. More precisely, there exists a stretch factor $\mu > 0$ such that almost surely
$$\f{T(x, y)}{|x - y|} \xra{|x - y| \tff} \mu,$$
where $T(x, y)$ denotes the smallest number of hops that are needed to connect $q(x)$ and $q(y)$, the points in $\CC_\ff^\la$ that are closest to $x$ and $y$ in Euclidean distance. 

Henceforth, we assume that $k \ua \ff$ and $\vli \da 0$ such that
\begin{align}
	\label{scale_la_eq}
	\vli (k/\mu)^d |B_1(o)|= \ns,
\end{align}
for some fixed $\ns > 0$, which may be interpreted as the expected number of sinks that are within $k$-hop range of a typical node at the origin. Indeed, since the stretch factor converts the $k$-hop distance to the Euclidean distance, any such sink is contained in the ball $B_{k/\mu}(o)$. Thus, the expected number of sinks is $\vli |B_{k/\mu}(o)|$.

To work out the impact of mobility cleanly, we need to take into account the sink-densities relation to the considered time horizon. More precisely, we investigate scalings of the form
\begin{align}
	\label{scale_t_eq}
	\vli(T) := T^{-\a},
\end{align}
for some parameter $\a > 0$ governing the sink density.

Concerning mobility, we assume that nodes choose random waypoints sequentially according to some isotropic probability measure $\k(\d v)$ and directly jump to them after exponentially distributed waiting times. In particular, $X(t)$ remains a Poisson point process with the same intensity. 
We assume that the trace of the coordinate-covariance matrix associated with random vectors from $\k(\d v)$ equals $d$. This normalization will later ensure convergence to a standard Brownian motion.

In Theorem \ref{thm_2}, we identify the connection measure after the scaling. To that end, let 
$\Xi^\ff := \{t\co X_0(t) {\lrsa_t} \ff\}$ be
the set of all times when the typical node is part of the infinite component of nodes. Similarly, 
$\Xi_o^\ff := \{t\co o {\lrsa_t} \ff\}$
denotes the set of all times when the static origin is part of the infinite component of nodes. Further we write $Y'(A)$ for the number of points of a point process $Y'$ in the measurable set $A\subset\R^d$.

\bet[Asymptotic weighted $k$-hop connection-interval measure]
\label{thm_2}
Let $\la > \la_c$ and assume the multi-scale regime encoded by \eqref{scale_la_eq} and \eqref{scale_t_eq}.\\[2ex]
{\bf Dense sinks.}
If $\a < d/2$, then, as $T \to \ff$,
  \begin{align}
	  \label{long_eq}
	  \tktb(\d \ell, \d t) \xrd \E[\de_{I_o(N)}(\d \ell)]\d t,
  \end{align}
 where  $I_o(N) := I(0, \Xi^\ff \cap (\cup_{j \le N}\Xi_o^{j, \ff}))$ with $N$ an independent Poisson random variable with intensity $\ns$ and $(\Xi_o^{j, \ff})_{j\ge 1}$ iid copies of $\Xi_o^\ff$.\\
{\bf Sparse sinks.}
If $\a > d/2$, then, as $T \to \ff$,
\begin{align}
	\label{short_eq}
	\tktb(\d \ell, \d t) \xrd \E[\de_{I_o(N)}(\d \ell) \ba N]\d t,
\end{align}
with all definitions as in the dense case. \\
{\bf Critical density.}
If $\a = d/2$, then, as $T \to \ff$,
  \begin{align}
	  \label{crit_eq}
	  \tktb \xrd \E[\de_{I_o(Y'(B_{\nsp}(W_t)))}(\d \ell) \ba Y'(B_{\nsp}(W_t))] \d t,
  \end{align}
where $\nsp :=(\ns/|B_1(o)|)^{1/d}$ and $Y'$ is a unit-intensity homogeneous Poisson point process and $W_t$ is a standard Brownian motion. 
\ent

We note that convergence of the distribution of random measures is defined by the convergence of integrals with respect to bounded continuous test functions. Hence, the statistics described in items (2) and (3) of Example \ref{stat_ex} should be truncated at some large value $M > 0$. In practice, when relying on such statistics as metrics for the network performance, the truncation at large $M$ is of little concern.

We conclude this section by expounding on how the asymptotic results presented in Theorem~\ref{thm_2} can be applied in the design and analysis of wireless networks. To that end, we concentrate on the dense case and present a simulation study where we describe the dependence of $\E[f(I_o(N))]$ on the expected number of in-range sinks $\ns$ for the three characteristics $f_1, f_2, f_3$ discussed in Example \ref{stat_ex}.

More precisely, the network nodes form a homogeneous Poisson point process with intensity $\la = 150$ in a $5 \ti 5$-sampling window, thus giving rise to a communication network with an expected number of 3,750 nodes and the communication radius is set to $0.1$. Although not needed for the computation of $\E[f(I_o(N))]$, for completeness, we note that the critical intensity for percolation is $\la_c \approx 143.7$ \cite{mertens} and that the stretch factor is $\mu \approx 8.1$ (own simulations). Observe that on the right-hand side in Theorem \ref{thm_2} the parameters $\vli,k,T$ are sent to the limit and therefore do not appear in the simulation. Figure \ref{net_fig} (left) shows a realization of this system. 

\bef
\includegraphics[width=\textwidth]{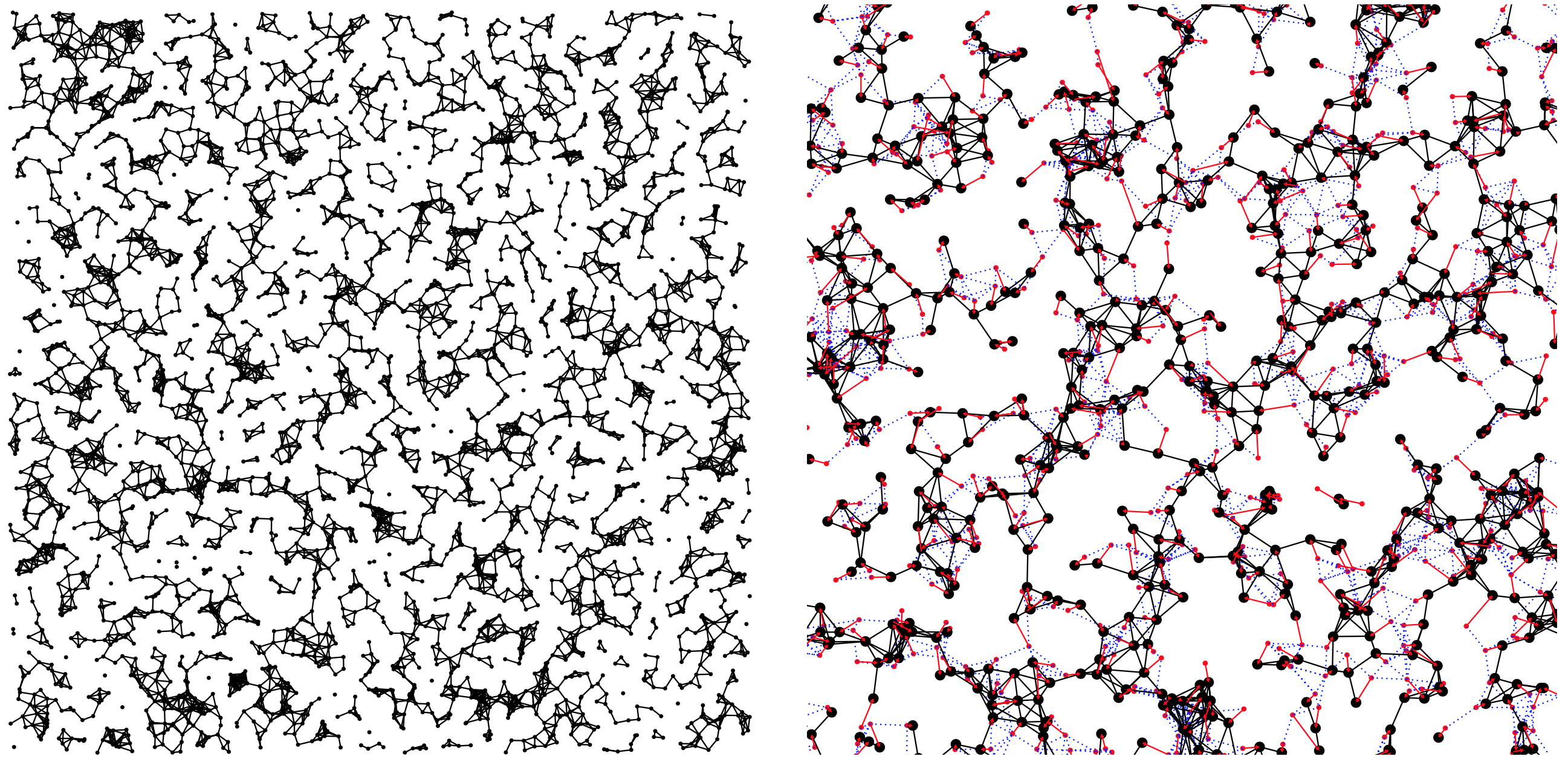}
 \caption{Sample of the simulated network at time 0 (left).  Node displacement after time 100 (right). In order to illustrate the mobility more clearly, we only show a cut-out of the full network. The network structure at time 0 is marked by blue dotted edges.}
 \label{net_fig}
\enf

Each node jumps according to a sequence of rate 1 exponential waiting times, and the jump locations are selected uniformly at random at distance $0.005$. Figure \ref{net_fig} (right) illustrates the changes of the network topology under these dynamics.

\bef
\includegraphics[width=.3\textwidth]{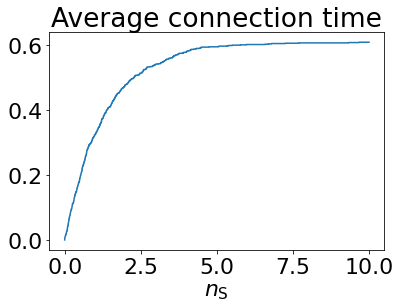}
\includegraphics[width=.358\textwidth]{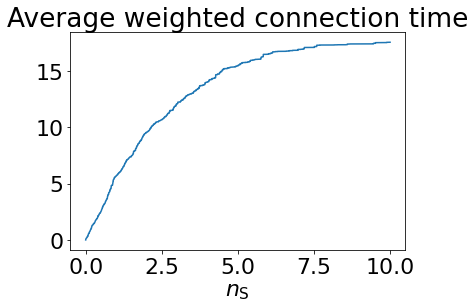}
\includegraphics[width=.325\textwidth]{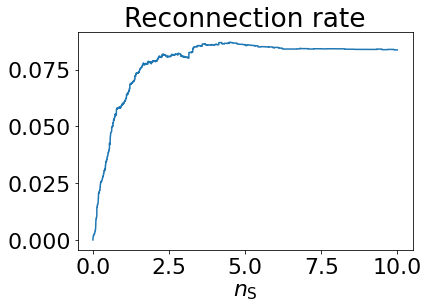}
\caption{Average connection time (left), average weighted connection time (center), and reconnection rate (right) based on 1,000 simulations.}
\label{sim_fig}
                       \enf

Now, we generate 1,000 realizations of this model and evaluate the quantity $\E[f(I_o(N))]$ appearing in Theorem \ref{thm_2} for the three different choices of the test function. Figure \ref{sim_fig} illustrates the Monte Carlo estimates for $\E[f_i(I_o(N))]$ as a function of $\E[N] = \ns$, the expected number of sinks in range.

Both $f_1$ and $f_2$ lead to increasing functions in $\ns$ approaching a saturation level as $\ns$ becomes large. We note that with respect to $f_1$ the system is already close to a saturation value of approximately $0.6$ for $\ns \ge 2$. That is, the typical node is connected approximately 60\% of the time. We note that this saturation value corresponds to the percolation probability at the considered intensity: if the typical node cannot form a connection to other devices outside a local neighborhood, then increasing the sink density does not help improving the total connectivity.

Concerning $f_2$, increasing $\ns$ may improve the system quality substantially, even beyond $\ns \ge2$. We hypothesize that this is due to the following effect. Although after $\ns \ge 2$, adding further sinks may not improve the total connection time substantially, still a few additional sinks may be enough to merge several smaller connection intervals into a very long one. This will boost the weighted connection time massively. 
 
	 Although the function $f_3(\ell, t)$ is neither increasing nor decreasing in $\ell$, Figure \ref{thm_2} still illustrates that it increases rapidly for small values of $\ns$ before reaching a plateau at around $0.075$ when $\ns \ge 2$. Loosely speaking, this means that if the typical node performs a jump on average once every second, then we see $0.075$ reconnections per second. Our interpretation is as follows: for low $\ns$, the number of reconnections is small simply because the typical device is in any event disconnected for most of the time. Then, as we increase $\ns$, the number of chances for the typical device to connect surges rapidly. Conceptually there is also a tendency in the opposite direction since an increase of $n_S$ eliminates some reconnections by merging two smaller intervals into a larger one. However, the plot suggests that the effect coming from the addition of new intervals is much stronger.

As a global conclusion of the considered simulation study, we may say that increasing $\ns$ past a very high value does not lead to further connectivity gains. However, we stress that the impact of increasing the number of admissible hops $k$ is a bit more subtle. Then,  also $\Xi^\infty$ becomes larger so that implementing this measure may boost connectivity even in situations where increasing $n_S$ does not help.

\section{Outline}
\label{out_sec}

The goal of Theorem \ref{thm_2} is to prove that the random measure $\tktb$ converges in distribution to a suitable random measure $\La$ on $[0, \ff) \ti [0, 1]$. As a preliminary step, we reduce this task to proving this convergence when integrating with respect to product-form test functions. More precisely, we consider integrals of the form
$$\tktb(gh) =  \int_0^1 g\big(I(tT, \Xi^k)\big) h(t) \d t$$
for continuous test functions $g\co \jtf$ and $h \co \itf$. The following auxiliary result then shows that it suffices to prove that the random variables $\tktb(gh)$ converge to $\La(gh)$ in distribution as $T\tff$.
%
%
\bel[Product-form test functions]
\label{prod_lem}
Let $\La$ be a random probability measure on $\off\ti [0, 1]$ such that $\tktb(gh)$ converges in distribution to $\La(gh)$ for every continuous $g\co \jtf$ and $h \co \itf$. Then, the random measure $\tktb$ converges weakly in distribution to the random measure $\La$.
\enl
We present the proof of Lemma~\ref{prod_lem} at the end of this section. 
%
The next preparatory step, is to show that we may replace the original expression $g\big(I(tT, \Xi^k) \big)$ by a finite-range approximation. More precisely, this involves three steps. First, we truncate long connection intervals. That is, we replace $I(tT, \Xi^k)$ by $I(tT, \Xi^k) \wedge M$ for a large truncation level $M > 0$. Second, we introduce $\de$-discretizations of $I(tT, \Xi^k)$, where we only check the connectivity at discretized time points. For this, we put
$$\idm(t, \Xi^k):= \sup_{\substack{a  \le b \in \jdm(t)\\ [a, b]\su \Xi^k}}(b - a),$$
where we rely on the finitely-many discrete time points 
$$\jdm(t):=  t + \{-\lc M/\de\rc \de, \dots, -\de, 0, \de, \dots, \lc M/\de\rc\de\}$$ around $t$.  Third, we replace the actual $k$-hop connection event by the event of percolation beyond bounded neighborhoods.  To that end, we write
$$\CC_L(t) := \big\{x \in \R^d\co x \lrsa_t\pa Q_L(x)\big\}$$
for the family of all locations percolating beyond an $L$-neighborhood via nodes in $X(t)$. 

First, in Section \ref{sec_finiterange}, we show that passing to the approximations incurs a negligible approximation error in the sense of Lemma \ref{approxi_lem} below. To make this precise, we let
$$Y(t, k) := Y(t, k, X_0) := B_{k/\mu}(X_0(t)) \cap Y$$
denote the family of all relevant sinks at time $t$. Then, we write $s\in\Xi^L(Y(t, k))$ if and only if $X_0(s) \in \CC_L(s)$ and $Y_j \in \CC_L(s)$ for some sink $Y_j \in Y(t, k)$. Moreover, for $L, M, \de > 0$, we set 
$$\idmm(tT) := \idm(tT, \Xi^L(Y(tT, k))).$$

\bepr[Finite-range approximation of connection intervals]
\label{approxi_lem}
Let $t \le 1$. Then,
$$\limsup_{M\tff}\limsup_{\de\da 0}\limsup_{L \tff}\limsup_{T\tff}\E\big[\big|I(tT, \Xi^k) - \idmm(tT)\big|\big] = 0.$$
\enpr


 After integrating with respect to the functions $g$ and $h$ and inserting the approximations, the first key step is to prove the following finite-range variant of Theorem \ref{thm_2}.

 \bepr[Approximate connection intervals]
\label{fr_thm}
Let $\la > \la_c$ and assume the multi-scale regime encoded by \eqref{scale_la_eq} and \eqref{scale_t_eq}.\\[2ex]
{\bf Dense sinks.}
If $\a < d/2$, then, for all $L,M,\de>0$, as $T \to \ff$,
  \begin{align}
          \label{long_apx_eq}
	  \int_0^1g\big(\idmm(tT) \big)h(t)\d t
	  \xrd \int_0^1\E[g\big(\ildm(N) \big)]h(t)\d t,
  \end{align}
 where  
 $$\ildm(N) := \idm\Big(0, \{s \in \R\co X_0(s) \in \CC_L(s)\} \cap \{s \in \R\co o \in \cup_{j \le N}\CC_L^{(j)}(s)\}\Big).$$
Here, $N$ is an independent Poisson random variable with intensity $\ns$ and the $\CC_L^{(j)}$ are iid copies of $\CC_L$.\\
{\bf Sparse sinks.}
If $\a > d/2$, then, for all $L,M,\de>0$, as $T \to \ff$,
\begin{align}
        \label{short_apx_eq}
	  \int_0^1g\big(\idmm(tT) \big)h(t)\d t \xrd \int_0^1\E\big[g\big(\ildm(N)\big) \ba N \big]h(t)\d t,
\end{align}
with all definitions as in the dense case. \\
{\bf Critical density.}
If $\a = d/2$, then, for all $L,M,\de>0$, as $T \to \ff$,
  \begin{align}
          \label{crit_apx_eq}
        	  \int_0^1g\big(\idmm(tT) \big)h(t)\d t   \xrd \int_0^1\E\big[g\big(\ildm(Y'(B_{\nsp}(W_t)))\big) \ba Y'(B_{\nsp}(W_t)) \big]h(t)\d t
  \end{align}
  where $\nsp :=(\ns/|B_1(o)|)^{1/d}$ and $Y'$ is a unit-intensity homogeneous Poisson point process and $W_t$ is a standard Brownian motion.
\enpr

We now explain how to deduce Theorem \ref{thm_2} from Propositions~\ref{approxi_lem} and~\ref{fr_thm}.

%
%
\bep[Proof of Theorem \ref{thm_2}]
We explain how to argue for $\a > d/2$, noting that the other two cases are similar. Let $F\co[0, 1] \rightarrow [0, 1]$ be Lipschitz with Lipschitz constant 1. We want to show that 
	  $$\lim_{T\tff }\E\Big[F\Big(\int_0^1g\big(I(tT, \Xi^k) \big)h(t)\d t \Big)\Big]=\E\Big[F\Big( \E[g\big(I_o(N) \big) \ba N]\int_0^1h(t)\d t\Big)\Big].$$
	  To that end, we consider the decomposition
	  \begin{align*}
		  &\Big|\E\Big[F\Big(\int_0^1g\big(I(tT, \Xi^k) \big)h(t)\d t \Big)\Big] -  \E\Big[F\Big( \E[g\big(I_o(N) \big) \ba N]\int_0^1h(t)\d t\Big)\Big]\Big|\\
		  &\le\Big|\E\Big[F\Big(\int_0^1g\big(\idm(tT, \Xi^L(Y(tT,k))) \big)h(t)\d t \Big)\Big] - \E\Big[F\Big(\E[g\big(\ildm(N) \big) \ba N] \int_0^1h(t)\d t\Big)\Big]\Big|   \\
		  &\phantom\le+\E\Big[\Big|F\Big(\int_0^1g\big(I(tT, \Xi^k) \big)h(t)\d t \Big) -  F\Big(\int_0^1g\big(\idm(tT, \Xi^L(Y(tT,k))) \big)h(t)\d t \Big)\Big|\Big]\\
		  &\phantom\le+\E\Big[\Big|F\Big( \E[g\big(\ildm(N) \big) \ba N]\int_0^1h(t)\d t\Big) - F\Big( \E[g\big(I_o(N) \big) \ba N]\int_0^1h(t)\d t\Big)\Big|\Big]\\
		  &\le\Big|\E\Big[F\Big(\int_0^1g\big(\idm(tT, \Xi^L(Y(tT,k))) \big)h(t)\d t \Big)\Big] - \E\Big[F\Big(\E[g\big(\ildm(N) \big) \ba N] \int_0^1h(t)\d t\Big)\Big]\Big|\\
		  &\phantom\le+\E\Big[\Big| g\big(I(0, \Xi^k) \big) -  g\big(\idm(0, \Xi^L(Y(0,k))) \big) \Big|\Big] +\E\Big[\Big|g\big(\ildm(N) \big)- g\big(I_o(N) \big)\Big|\Big],
	  \end{align*}
  where in the last inequality we used that the model is time-stationary. By Proposition~\ref{fr_thm}, the first expression on the right-hand side tends to 0 as $T \tff$. Since $g$ is continuous, Proposition \ref{approxi_lem} allows us to choose $M, L, \de > 0$ such that the second expression becomes arbitrarily small as $T \tff$. Finally, a similar argument applies to the third contribution, thereby concluding the proof.
\enp

The key step in the proof of Proposition \ref{fr_thm} is the conditional second-moment method. Hence, we need to control conditional expectations and covariances for expressions like $g\big(\idmm(tT)\big)$, $t \le 1$.

%
%
\bepr[Asymptotic conditional decorrelation]
        \label{decorr_pr}
        Let $0 < t < 1$. Then, under the scalings \eqref{scale_la_eq} and \eqref{scale_t_eq},
\begin{align}
        \label{decorr_eq}
\lim_{T\to\ff}\E\big[\Cov\big(g\big(\idmm(0)\big), g\big(\idmm(tT)\big) \ba X_0, Y\big)\big]= 0.
\end{align}
\enpr

%
%
To describe concisely the asymptotic conditional expectation given $X_0$ and $Y$, we need a more explicit representation of $\idmm$. By definition, $\idmm(tT)$, is determined by certain finite-range percolation events of the typical node $X_0$, and of sinks in range $Y(tT, k)$  at the times $\jdm(tT)$. More precisely, to make this dependence explicit, we will also write 
$$\idm\Big(\{s \in \jdm(tT)\co X_0(s) \in \CC_L(s)\}, \{s \in \jdm(tT)\co Y_j \in \CC_L(s)\}_{{Y_j \in Y(tT, k)}} \Big)$$
instead of $\idmm(tT)$.


%
%
\bepr[Asymptotic conditional expectation]
        \label{xyc_pr}
        Let $ t \le 1$. Then, under the scalings \eqref{scale_la_eq} and \eqref{scale_t_eq}, we have the $L^1$-convergence
\begin{align*}
	\E[g(\idmm(tT)) \ba X_0, Y]
	\xra{T\tff}  S\big(\{X_0(s) - X_0(tT)\}_{s \in \jdm(tT)}, \#Y(tT, k)\big),
\end{align*}
where
\begin{align}
	\label{sp_eq}
S\big(\{x_s\}_{s \in \jdm(tT)}, n\big) := 
	\E\big[g\big(\idm\big(\{s \co x_s \in \CC_L(s)\}, \{s\co o \in \CC_L^{(j)}(s)\}_{j \le n} \big\}\big)\big)\big]
\end{align}
\enpr
We prove Proposition \ref{decorr_pr} and \ref{xyc_pr} in Section \ref{decorr_sec}. After that, we will elucidate in Sections~\ref{dense_sec}--\ref{crit_sec} below how to complete the proof of Proposition~\ref{fr_thm} in the different regimes for the parameter $\a$.
We conclude this section with a proof of Lemma \ref{prod_lem}.

\bep[Proof of Lemma \ref{prod_lem}]
We follow the proof of \cite[Theorem 4.11]{kallenberg}. First, we show the tightness of the random probability measures $\tktb$. To that end, we note that by assumption, the random variable $\tktb(gh)$ converges to $\La(gh)$ in distribution for every $g, h$ as above. Thus, also $\lim_{T \tff}\E[\tktb(gh)] = \E[\La(gh)]$ for every $g, h$ as above. Since the underlying space $\off \ti [0, 1]$ is a product space, we deduce from \cite[Lemma 4.1]{kallenberg} that the expectation measure $\E[\tktb]$ converges weakly to $\E[\La]$. Thus, by Prohorov's theorem (see \cite[Theorem 4.2]{kallenberg}), $\lim_{M \tff}\sup_{T \ge 1}\E[\tktb([M, \ff) \ti [0, 1])] = 0$. By \cite[Theorem 4.10]{kallenberg} this yields tightness of the random measures $\{\tktb\}_{T \ge 1}$. Now, we conclude as in the proof of \cite[Theorem 4.11]{kallenberg}.
\enp

\section{Proofs}\label{Sec_Pro}

In Section, \ref{sec_finiterange}, we establish the finite-range approximation of the connection intervals from Proposition \ref{approxi_lem}.
In Section \ref{decorr_sec}, we establish the asymptotic conditional?expectations and covariances statements of the Propositions \ref{decorr_pr} and \ref{xyc_pr}.
\subsection{Finite-range approximation -- proof of Proposition \ref{approxi_lem}}
\label{sec_finiterange}

First, we may truncate long connection intervals. More precisely, the typical connection intervals are tight.
%
%
\bel[Tightness]
\label{tight_lem}
Under  \eqref{scale_la_eq} and \eqref{scale_t_eq}, the random variables $\{I(0, \Xi^k)\}_{k \ge 1}$ are tight.
\enl

 Second, the discretization error vanishes as $\de \da 0$.

%
%
\bel[Discretization]
\label{disc_lem}
Let $M > 0$. Then, under  \eqref{scale_la_eq} and \eqref{scale_t_eq},
$$\lim_{\de \da 0}\limsup_{k \to \ff}\E\big[|\idm(0, \Xi^k) - I(0, \Xi^k) \wedge M|\big] = 0.$$
\enl

Finally, we show how to approximate the $k$-connection event by the percolation outside finite boxes. 

%
%
\bel[Finite-range percolation]
\label{fp_lem}
Let $M, \de  > 0 $. Then, under  \eqref{scale_la_eq}, \eqref{scale_t_eq},
$$\lim_{L \tff}\limsup_{k \to \ff}\E\big[|\idm(0, \Xi^k) - \idmm(0)|\big] = 0.$$
\enl

%
%
In the rest of this section, we prove Lemmas \ref{tight_lem}--\ref{fp_lem}. We start with Lemma \ref{fp_lem} since it follows from a short argument based on the shape theorem for continuum first-passage percolation. The latter allows to replace the $k$-hop connection event by suitable percolation events. 
\bep[Proof of Lemma \ref{fp_lem}]
First, note that $\lim_{k \tff} \P\big(\text{$Y(t', k) = Y(0, k)$ for all $t' \in [ - M,  M]$}\big)  = 1$. Thus, by the shape theorem for continuum first-passage percolation in the form of \cite[Lemma 15]{mscale}, 
$$\lim_{k\tff} \E\big[\big|\idm(0, \Xi^k) - \idm(0, \Xi^*(X_0, Y(0, k) ) )\big|\big] = 0, $$
where $t'\in \Xi^*(X_0, Y(0, k))$ if $X_0(t') \in \CC_\ff(t')$ and $Y_j \in \CC_\ff(t')$ for some sink $Y_j \in Y(tT, k)$.

Now, by the uniqueness of the connected component, $X_0(t') \in \CC_\ff(t')$ if and only if $X_0(t') \in \CC_L(t')$ for every $L \ge 1$ and similarly $Y_j \in \CC_\ff(t')$ if $Y_j \in \CC_L(t')$ for every $L \ge 1$. Hence, we may replace $\Xi^*(X_0, Y(0, k))$ by $\Xi^L(Y(0, k))$, thereby concluding the proof.
\enp

Next, we prove Lemma \ref{tight_lem}. As a pivotal observation, we note that since nodes perform a random walk, the node movement is diffuse in the sense that it is highly unlikely that after a time of order $T$, the typical node is contained in a set of bounded diameter.

%
%
\bel[Diffuseness of node locations]
\label{diff_node_lem}
Let $t, L > 0$. Then, under  \eqref{scale_la_eq} and \eqref{scale_t_eq}
$$\lim_{T\tff}\sup_{x \in \R^d} \P\big(X_0(tT) \in Q_L(x)\big) = 0.$$
\enl
\bep
Fix $\e > 0$ and note that we may assume $L > 1$. First, by the central limit theorem, $X_0(tT)/\sqrt T$ converges in distribution to a Gaussian vector $Z$. Then, for $M \ge 1$ partition the box $Q_M$ into $M^{2d}$ boxes $Q_{1/M}(z_{1, M}), \dots, Q_{1/M}(z_{M^{2d}, M})$ of side lengths $1/M$ centered at points $z_{1, M}, \dots, z_{M^{2d}, M} \in Q_M$. In particular, we may fix $M$ so large that $\P(Z \not \in Q_M) < \e$ and $\sup_{i \le M^{2d}}\P(Z \in Q_{3/M}(z_{i, M})) < \e$.

Now, if $T$ is so large that $\sqrt T/M \ge L$, then every Borel set $B$ of diameter at most $L$ is contained in $\R^d\sm Q_{\sqrt TM}$ or in $Q_{3\sqrt T/M}(z_{i, M})$ for some $i \le M^{2d}$. Hence, invoking the distributional convergence of $X_0(tT) / \sqrt T$ concludes the proof.
\enp

The key observation for the tightness assertion in Lemma \ref{tight_lem} is to use that a large connection interval means that the typical node needs to be in the unbounded connected component for many temporally distant time steps.

%
%
\bep[Proof of Lemma \ref{tight_lem}]
We show tightness of the connection interval $I^+(0, \Xi^k) := I(0, \Xi^k \cap \off)$ for positive connection times. The arguments for negative connection times are symmetric.

Let $\e > 0$. Then, by discretization, for any $n_0, t_0 \ge 1$,
$$\P\big(I^+(0, \Xi^k) > n_0t_0\big) \le \P\big(\{0, t_0, \dots, n_0t_0\} \su  \Xi^k \big).$$
Since $\lim_{L \tff} \th_L :=\lim_{L \tff}  \P(o \in \CC_L) = \th < 1$, we may fix $L_0, n_0 \ge 1$ such that $\th_{L_0}^{n_0} < \e$.

Next, if $nt_0 \in \Xi^k$ and there are no sinks in an $L_0$-neighborhood around $X_0(nt_0)$, i.e., if $Y \cap Q_{L_0}(X_0(nt_0)) = \es$, then it is possible to percolate beyond that neighborhood, i.e., $X_0(nt_0)\in \CC_{L_0}(nt_0)$. Thus, it suffices to produce $t_0 \ge 1$ such that
$$\P\big(X_0(nt_0)\in \CC_{L_0}(nt_0) \text{ for all $n \le n_0$}\big) \le 2\e.$$
To that end, we let
$$X^{-,nt_0} := \big\{X_i \in X\co \{n \} = \{n' \le n_0\co X_i(n't_0)\in Q_{L_0}(X_0(n't_0))\}\big\}$$
denote the family of nodes that are contained in the $L_0$-neighborhood of the typical node at time $nt_0$ but not at any other discretized times. We also introduce $\CC_L^-(t)$ in the same way as $\CC_L(t)$ with the only difference that percolation exclusively relies on nodes in $X^{-,t}$.
In particular, by the independence property of the Poisson point process, when conditioned on $X_0$,
$$\P\Big(X_0(nt_0)\in \CC_{L_0}^-(nt_0)  \text{ for all $n \le n_0$} \ba X_0\Big) = \prod_{n \le n_0} \P\big(X_0(nt_0) \in \CC_{L_0}^-(nt_0)  \ba X_0\big).$$
Since  $\P\big(X_0(nt_0) \in \CC_{L_0}(nt_0) \ba X_0 \big)=\th_{L_0}$, it therefore suffices to show that almost surely,
$$\lim_{t_0 \tff}\P\big(X_0(nt_0) \in \CC_{L_0}(nt_0)\ba X_0  \big) - \P\big(X_0(nt_0) \in \CC_{L_0}^-(nt_0)  \big) = 0.$$
Now, the difference on the left-hand side is bounded above by
\begin{align*}
	&\P\big(X(nt_0)\cap Q_{L_0}(X_0(nt_0)) \ne X^-(nt_0)\cap Q_{L_0}(X_0(nt_0))\big)\\
	&\quad\le \E\big[\big|(X(nt_0) \sm X^-(nt_0))\cap Q_{L_0}(X_0(nt_0))\big|\big],
\end{align*}
which tends to 0 as $t_0 \tff$ by Lemma \ref{diff_node_lem}.
\enp

The main idea in the proof of Lemma \ref{disc_lem} is that for small $\de$, only very few nodes move within the time interval $[i\de, (i + 1)\de]$, and therefore we do not need to rely on them when establishing $k$-hop connections.

In order to make this more precise, note that by the thinning theorem, the intensity of nodes that are moving in an interval of the form $[i_0\de_0, (i_0 + 1)\de_0]$ form a Poisson point process $X^{i_0, \de_0}$ with intensity $(1 - e^{-\de_0})\la$, see \cite[Theorem 5.8]{poisBook}. In particular, for sufficiently small $\de_0$, this process is still in the super-critical phase and we let $\ccn(i_0\de_0)$ denote the associated  unique unbounded connected component for continuum percolation. We write $E_{k, \de_0}^{\ms{glob}}$ for the event that all pairs of nodes $X_j \in \ccn(i_0\de_0) \cap B_{k^{1/(2d)}}$ and $X_{j'} \in \ccn(i_0\de_0) \cap B_{(1 - \de_0)k/\mu}$ are connected in at most $k - \sqrt k$ hops for every $i_0 \in \Z$ with $|i_0| \le M/\de_0$. The first step in the proof of Lemma \ref{disc_lem} is to show that these global connections occur with a high probability provided that $\de_0$ is sufficiently small.

\bel[Global connections]
\label{glob_lem}
If $\de_0 > 0$ is sufficiently small, then
$\lim_{k \tff}\P(E_{k, \de_0}^{\ms{glob}}) = 1.$
\enl

In order to complete the global paths to a full $k$-hop connection, we still require short local paths leading up to the unbounded connected components. To that end, we will need to refine the above $\de_0$-discretization and will rely on variants of $\Xi^k$ that are locally determined.  More precisely, for $L > 0$ we write that $t \in \Xi^{L, \de_0}$ if and only if $X_0$ connects at time $t$ to some node of $\ccn(\lf t/\de_0\rf\de_0) \cap Q_L(X_0(t))$ in at most $L^{2d}$ hops. Then, an essential step is to show that these local paths connect to the global path with high probability. To that end, we couple $L$ to the size of a finer discretization by putting $L_\de := \de^{-1/(2d)}$. Now, we define
$$E_{i, \de}^{\ms{loc}}:= \Big\{X_0(i\de)\in \CC_{L_\de}(i\de),\, X_0((i + 1)\de)\in \CC_{L_\de}((i + 1)\de),\, [i\de, (i + 1)\de]\not\su \Xi^{L_\de, \de_0}\Big\}$$
as  the event that there are both $X_0(i\de)$ and $X_0((i + 1)\de)$ percolate beyond an $L_\de$-neighborhood but that $[i\de, (i + 1)\de] \not \su \Xi^{L_\de, \de_0}$.

\bel[Local connections]
\label{loc_lem}
Let $M > 0$. Then,
$\sup_{|i| \le M/\de}  \P\big(E_{i, \de}^{\ms{loc}}\big) \in o(\de).$
\enl

%
%
\bep[Proof of Lemma \ref{disc_lem}]
The task is to show that
\begin{align}
        \label{disc_asL_eq}
        \lim_{\de \da0}\limsup_{k\to\ff}\P\Big(\bigcup_{{|i| \le M/\de}}E_{i, k,\de}\Big) = 0,
\end{align}
where
$$E_{i,k, \de}: = \big\{\{i\de, (i + 1)\de\} \su \Xi^k, [i\de, (i + 1)\de]\not\su \Xi^k\big\}$$
is the event of having $k$-connections at $i\de$ and $(i + 1)\de$ but not within the entire interval $[i\de, (i + 1)\de]$.
Note that if $i\de \in \Xi^k$, then $X_0(i\de) \in \CC_{L_\de}(i\de)$ unless some sink lies in $Q_{L_\de}(X_0(i\de))$. However, since the sink intensity tends to 0 as $k \tff$ so does the probability of the latter event. Moreover, by Lemma \ref{loc_lem},
$$\limsup_{\de\da0}\P\big(\bigcup_{|i| \le M/\de}E_{i, \de}^{\ms{loc}}\big) =2M\limsup_{\de\da0} \tf1\de\sup_{|i| \le M/\de}\P\big(E_{i, \de}^{\ms{loc}}\big) = 0.$$
Noting that a similar argument applies when replacing $X_0(i\de)$ by one of the sink nodes in $B_{k/\mu}(X_0(i\de))$ thus concludes the proof of identity \eqref{disc_asL_eq}.
\enp

Finally, we prove Lemmas \ref{glob_lem} and \ref{loc_lem}.

%
%
\bep[Proof of Lemma \ref{glob_lem}]
The key ingredient in the proof is the continuity of the stretch factor with respect to the intensity $\la$ of the underlying Poisson point process. While in first-passage percolation on the lattice, results in this vein are classical and hold under very general conditions \cite{stretch}, the question in the continuum only requires a small adaptation. We only show that $\limsup_{\la' \to \la}\mu_{\la'} \le \mu_{\la}$ since this will be sufficient for the proof of Lemma~\ref{glob_lem}. Once this assertion is established, we invoke the shape theorem for continuum first-passage percolation in the form of  \cite[Lemma 15]{mscale} in order to deduce that $\lim_{k \to \ff} \P(E_{k, \de_0}^{\ms{glob}}) = 1$ for sufficiently small $\de_0$.

To that end, set $\mu_{\la, n} := n^{-1}\E[T_\la(q_{ \la}(o), q_{ \la}(ne_1))]$, where $q_{ \la}(x)$ is the point in $\CC_\ff^\la$ that is at smallest Euclidean distance to $x \in \R^d$ and $T_\la(x,y)$ denotes the graph distance between points $x,y\in \CC_\ff^\la$. Then, $T_{\la'}(q_{ \la'}(o), q_{ \la'}(ne_1))$ converges almost surely to $T_\la(q_{ \la}(o), q_{ \la}(ne_1))$ as $\la'\to\la$. Hence fixing some super-critical $\la_0$, it suffices to show that the path lengths $T_{\la'}(q_{ \la'}(o), q_{ \la'}(ne_1))$ are uniformly integrable as $\la' \in [\la_0, \la]$.

To achieve this goal, we note that
$$T_{\la'}(q_{ \la'}(o), q_{ \la'}(ne_1)) \le T_{\la_0}\big(q_{ \la_0}(o), q_{ \la_0}(ne_1)\big) + T_{\la'}\big(q_{ \la'}(o), q_{ \la_0}(o)\big) + T_{\la'}\big(q_{ \la'}(ne_1), q_{ \la_0}(ne_1)\big).$$
Hence, by stationarity, it suffices to establish the uniform integrability of $T_{\la'}\big(q_{ \la'}(o), q_{ \la_0}(o)\big)$  for $\la' \in [\la_0, \la]$. 
For $M > 0$, say that the box $Q_M$ is \emph{$M$-good} if $\CC_\ff^{\la_0} \cap Q_{M/2} \ne \es$ and if the unique component of the Gilbert graph on $X^\la$ in $Q_M$ of diameter more than $M/8$ contains  $\CC_\ff^{\la_0} \cap Q_M$.
In particular, if $Q_M$ is $M$-good, then $q_{ \la'}(o)$ and $q_{ \la_0}(o)$ are connected by a path in $Q_M$ so that $T_{\la'}\big(q_{ \la'}(o), q_{ \la_0}(o)\big) \le X^\la(Q_M)$. Thus, for any $a > 0$,
$$\P\big(T_{\la'}\big(q_{ \la'}(o), q_{ \la_0}(o)\big) > a\big) \le \P(X^\la(Q_{a^{1/(2d)}}) > a ) + \P(\text{$Q_{a^{1/(2d)}}$ is not $a^{1/(2d)}$-good}).$$
Now, we conclude from the quantitative uniqueness in the form of \cite[Theorem 2]{uniqFin} that $Q_M$ is $M$-good with high probability. Hence, the right-hand side becomes arbitrarily small for large $a$, thereby concluding the proof of the uniform integrability.
\enp

The proof of Lemma \ref{loc_lem} relies heavily on the observation that it is highly unlikely to see a time interval of a small length $\de$ where two or more nodes are moving. To make this precise, we let
$$E_{i, \de}^{\ms{disc}} := \big\{X(t) \cap Q_{L_\de}(X_0(t))\in \big\{X(i\de) \cap Q_{L_\de}(X_0(i\de))\co j \in \{i, i +1\}\big\} \text{  $\forall t \in [i\de, (i + 1)\de]$}\big\}$$
be the event that for any $t \in [i\de, (i + 1)\de]$ the configuration of $X(t) \cap Q_L(X_0(t))$ coincides with $X(i \de) \cap Q_{L_\de}(X_0(i\de))$ or $X((i + 1)\de) \cap Q_{L_\de}(X_0((i + 1)\de))$. By stationarity, $\P\big((E_{i, \de}^{\ms{disc}})^c\big)$ does not depend on $i$.

\bel[$E_{i, \de}^{\ms{disc}}$ occurs whp]
\label{discc_lem}
Let $M > 0$. Then,
$\sup_{|i| \le M/\de}  \P\big((E_{i, \de}^{\ms{disc}})^c\big) \in o(\de).$
\enl
\bep
First, the probability that the typical node $X_0$ jumps at least twice in $[i\de, (i + 1)\de]$ is $1 - e^{-\de} - \de e^{-\de} \in O(\de^2)$. Hence, we may assume that $X_0(t) \in \{X_0(i\de), X_0((i + 1)\de)\}$ for all $t \in [i\de, (i + 1)\de]$. We now distinguish between the cases whether or not the typical node moves in the interval $[i\de, (i + 1)\de]$. 
\smallskip

\noindent\boldsymbol{$X_0(i\de) \ne X_0((i + 1)\de).$}
First, note that $\P\big(X_0(i\de) \ne X_0((i + 1)\de)\big) = 1- e^{-\de} \in O(\de)$. Hence, it suffices to show that conditioned on $X_0$, the probability that either of $X(t)\cap Q_{L_\de}(X_0(i\de))$ or $X(t)\cap Q_{L_\de}(X_0((i + 1)\de))$ changes within the time interval $t \in [i\de, (i + 1)\de]$ is of order $o(1)$ as $\de \da0$. By time reversibility, it suffices to consider the case $i\de$;  by independence of $X_0$ and $X$, we may assume that $X_0(i\de) = o$.

The expected number of nodes of $X(i\de) \cap Q_{L_\de}(o)$ moving in the time interval $[i\de, (i + 1)\de]$ is at most $L_\de^d(1 - e^{-\de})$, and therefore of order $o(1)$. In order to bound the number of nodes entering $Q_{L_\de}(o)$ from the outside, we apply the mass-transport principle \cite{mtp1,mtp2}. More precisely, for $z, z'\in \Z^d$, we let $S(z, z')$ denote the total number of visits in $  Q_{L_\de}(L_\de z')$ within the time interval $[i\de, (i + 1)\de]$ of nodes that are contained in $ Q_{L_\de}(L_\de z )$ at time $i\de$. Then, the number of entering nodes is bounded above by the incoming mass at the origin. Conversely, the expected outgoing mass at the node $o$ is at most $L_\de^d\de \in o(1)$. Since the mass-transport principle  implies that the expected incoming mass equals the expected outgoing mass, we conclude the proof in the case $X_0(i\de) \ne X_0((i + 1)\de)$. 

\noindent\boldsymbol{$X_0(i\de) = X_0((i + 1)\de).$} In this case, we need to show that the probability that there is more than one change of $X(t)\cap Q_{L_\de}(X_0(i\de))$ in the time interval $[i\de, (i + 1)\de]$ is of order $o(\de)$. If there is more than one change, then this may be either because several nodes move or because some node moves multiple times.

First, consider the situation involving multiple nodes. Similarly, as in the previous case, conditioned on $X_0$, the collection of nodes in  $X(i\de) \cap Q_{L_\de}(X_0(i\de))$ moving in $[i\de, (i + 1)\de]$ is a Poisson point process with intensity $1 - e^{-\de}$. Hence, the probability of more than one of those nodes moves is of order $O(\de^2)$. Next, conditioned on $X_0$, the nodes entering $Q_{L_\de}(X_0(i\de))$ form a Poisson point process. Again, using the mass-transport principle the intensity of this process is of order at most $O(L_\de^d\de)$. Hence, the probability of seeing more than one entering node is of order at most $O(L_\de^{2d}\de^2)$ and therefore in $o(\de)$. Finally, it may happen that in the interval $[i\de, (i + 1)\de]$ at least one node from $X(i\de) \cap Q_{L_\de}(X_0(i\de))$ moves, and at least one node from $X(i\de) \sm Q_{L_\de}(X_0(i\de))$ enters $Q_{L_\de}(X_0(i\de))$. Conditioned on $X_0$, these two processes are independent, and by the arguments derived above there intensities are of order at most $O\big(L_\de^d (1 - e^{-\de})\big)$ and $O\big(L_\de^d\de\big)$, respectively. Thus, the asserted probability is of order at most $O\big(L_\de^{2d}(1 - e^{-\de})\de\big)$, as asserted.

Second, consider the situation of some nodes moving multiple times. Here, we again apply the mass-transport principle similarly to the setting where $X_0(i\de) \ne X_0((i + 1)\de)$. More precisely, for $z, z'\in \Z^d$, we now let  $S'(z, z')$ denote the total number of visits in $  Q_{L_\de}(L_\de z')$ within the time interval $[i\de, (i + 1)\de]$ of nodes that are contained in $ Q_{L_\de}(L_\de z )$ at time $i\de$ and that jump at least twice in the interval $[i\de, (i + 1)\de]$. Then, the number of relevant nodes moving multiple times is bounded above by the incoming mass at the origin. Conversely, the expected outgoing mass at the node $o$ is of the order at most  $O(L_\de^{2d}\de^2) \in o(\de)$. Hence, another application of the mass-transport principle concludes the proof.
\enp

%
%
\bep[Proof of Lemma \ref{loc_lem}]
For $\de > 0$ with $\de_0/\de \in \Z$ and $|i| \le M/\de$ let $i_0 = i_0(i, \de)$ be such that $[i\de, (i + 1)\de] \su [i_0\de_0, (i_0 + 1)\de_0]$.

A key ingredient is the strong quantitative uniqueness in the form of \cite[Theorem 2]{uniqFin}: For $\de > 0$  and $i \in \Z$ write $E_{i, \de}^{\ms{uniq}}$ for the event that for every $j \in \{i, i + 1\}$, it holds that
\been
\im inside $Q_{L_\de}(X_0(j\de))$ there exists a unique component of nodes in $X(j\de)$ of diameter at least $L_\de/8$;
\im inside $Q_{L_\de}(X_0(j\de))$ there exists a unique component of nodes in $X^{i_0, \de_0}(j\de)$ of diameter at least $L_\de/8$. Moreover, this component intersects $\ccn(i_0\de_0)$.
\enen
Then, \cite[Theorem 2]{uniqFin} provides a $c > 0$ such that $ 1 - \P\big(E_{i, \de}^{\ms{uniq}}\big) \le \exp(-cL_\de)$. Again by stationarity, the latter probability does not depend on $i$.

We now claim that  $E_{i, \de}^{\ms{loc}}$ cannot occur under the event 
$$\big\{X^\la(Q_{L_\de}(X_0(i\de))) \le L_\de^{2d}\big\} \cap \big\{X^\la(Q_{L_\de}(X_0((i + 1)\de))) \le L_\de^{2d}\big\} \cap E_{i, \de}^{\ms{uniq}} \cap E_{i, \de}^{\ms{disc}},$$ 
which will conclude the proof of the lemma since each of the probabilities $1 - \P(E_{i, \de}^{\ms{disc}})$, $1 - \P(E_{i, \de}^{\ms{uniq}})$ and $1 - \P(X^\la(Q_{L_\de}(X_0(i\de))) \le L_\de^{2d})$ is at most $o(\de)$.
Now, suppose that $X_0(i\de) \in \CC_{L_\de}(i\de)$, $X_0((i + 1)\de) \in \CC_{L_\de}((i + 1)\de)$  and let $t \in [i\de, (i + 1)\de]$ arbitrary.  
In particular, under the event $E_{i, \de}^{\ms{uniq}}$, we conclude that $X_0(j\de)$ connects within $Q_{L_\de}(X_0(j\de))$ to a vertex in $\ccn(i_0)$ for every $j \in \{i, i + 1\}$. Then, $i\de, (i + 1)\de \in \Xi^{L_\de, \de_0}$. Again, under the event $E_{i, \de}^{\ms{disc}}$ this means that $t \in \Xi^{L_\de, \de_0}$. Hence, the event $E_{i, \de}^{\ms{loc}}$ does not occur.
\enp

\subsection{Asymptotic conditional expectations and covariances -- proofs of Propositions \ref{decorr_pr} and \ref{xyc_pr}}
\label{decorr_sec}

The key to achieving asymptotic independence will be the diffuseness of the node movement. More precisely, it is highly unlikely to find a node contained in two specific neighborhoods at two distant points in time.

%
%
\bel[Asymptotic decorrelation of node locations]
\label{diff_lem}
Let $\de, M, K, t > 0$. Then,
$$\lim_{T\tff}\sup_h \sup_{\vp, \vp'}\Cov\Big(h\big(\{X(s) \cap Q_K(\vp)\}_{s \in \jdm(0)  }\big), h\big(\{X(s) \cap Q_K(\vp')\}_{s \in \jdm(tT)}\big)\Big) = 0,$$
where the suprema run over all measurable $[0, 1]$-valued functions $h$ and all  $K$-element subsets $\vp, \vp' \su \R^d$.
\enl

Before establishing Lemma \ref{diff_lem}, we elucidate how it enters the proof of Proposition \ref{decorr_pr}.

%
%
\bep[Proof of Proposition \ref{decorr_pr}]
First, by choosing $K\ge 1$ sufficiently large, we may assume that both $Y(0, k)$ and $Y(tT, k)$ have at most $K$ elements. More precisely, we replace 
$$\E\big[\Cov\big(g\big(\idmm(0)\big), g\big(\idmm(tT)\big) \ba X_0, Y\big)\big]$$
by 
$$\E\big[\Cov\big(g\big(\idmm(0)\big), g\big(\idmm(tT)\big) \ba X_0, Y\big) \one\{\#Y(0, k) \vee \#Y(tT, k) \le K\}\big].$$

After constraining the number of elements of $Y(0, k)$ and $Y(tT, k)$ they become eligible choices for $\vp$ and $\vp'$ in Lemma \ref{diff_lem}. Thus, that result allows us to conclude the proof.
\enp

Next, we establish the asymptotic representation of the conditional expectation from Proposition \ref{xyc_pr}. The proof idea is to use that the sinks are so sparse that no moving node can visit neighborhoods of two distinct sinks.

%
%
\bep[Proof of Proposition \ref{xyc_pr}]
First, by stationarity, we may assume that $t = 0$. Next, we let
$$E_k^{\ms{inert}} := \big\{ \sup_{s \in \jdm(0)}|X_0(s)| \le \sqrt k /4\big\}$$
denote the high-probability event that the typical node does not move further than $\sqrt k$ within the times $s \in \jdm(0)$. Moreover, we note that the high-probability event
$$E_k^{\ms{disj}} := \big\{|x - x'| \ge \sqrt k \text{ for all $x \ne x' \in Y(0, k) \cup \{o\}$}\big\}$$
implies the disjointness of the $L$-neighborhoods relevant for the percolation events encoded in $s \in \Xi^L(Y(tT, k))$. Also the events $E_k^{\ms{disj}}$ occur whp.

For $s \in \jdm(0)$, we now let 
$$X^{\ms{exc}} := \big\{X_i \in X \co \sup_{s' \in \jdm(0)}|X_i(s') - X_i(0)| \ge \sqrt k /4\big\}$$
denote the family of all nodes that move a distance further than $\sqrt k /4$ within the times. By the thinning theorem \cite[Corollary 5.9]{poisBook}, $X^{\ms{exc}}$ is a homogeneous Poisson point process with a vanishing intensity as $k \tff$. 

Next, we define the percolation sets $\CC_L^-(s)$ precisely as $\CC_L(s)$ except that the connections are only formed through nodes in $X \sm X^{\ms{exc}}$ instead of $X$. Thus, by relying on these modified percolation sets, we can introduce the sets $\Xi^{-, L}(Y(0, k))$.

Now, when conditioning on $X_0$ and $Y$, then under the event $E_k^{\ms{disj}} \cap E_k^{\ms{inert}}$, by the independence property of the Poisson point process on disjoint sets, the random vectors $Z_0 := \big\{\one\{X_0(s) \in \CC_L^-(s)\}\big\}_{s \in \jdm(0)}$ and $Z(Y_i) :=\big\{\one\{ Y_i \in \CC_L^-(s)\}\big\}_{s \in \jdm(0)}$, $Y_i \in Y(0, k)$ are independent. Therefore, under $E_k^{\ms{disj}} \cap E_k^{\ms{inert}}$ we have 
$$\E[g(\idm(0, \Xi^{-, L}(Y(0, k)))) \ba X_0, Y] =   S^-\big(\{X_0(s) \}_{s \in \jdm(0)}, \#Y(0, k)\big),$$
where 
$$S^-\big(\{x_s\}_{s \in \jdm(0)}, n\big) :=
\E\big[g\big(\idm\big(\{s \co x_s \in \CC_L^-(s)\}, \{s\co o \in \CC_L^{-,(j)}(s)\}_{j \le n} \big\}\big)\big)\big].
$$
Since the events $E_k^{\ms{disj}} \cap E_k^{\ms{inert}}$ occur whp, and since the intensity of $X^{\ms{exc}}$ vanishes as $k \tff$, we deduce that we may replace $S^-$ by $S$, thereby concluding the proof.
\enp

Lemma \ref{diff_node_lem} is the key ingredient for the proof of Lemma \ref{diff_lem}.

%
%
\bep[Proof of Lemma \ref{diff_lem}]
By the thinning theorem, the family
$$X^{t, \vp} := \big\{X_i \in X\co X_i(s) \in Q_K(\vp) \text{ for some $s \in \jdm(t)$}\big\}$$
of nodes contained in $Q_K(\vp)$ at some time point $s \in \jdm(t)$ forms a Poisson point process. Moreover,
\begin{align*}
        &\Cov\Big(h\big(\{X(s) \cap Q_K(\vp)\}_{s \in \jdm(0)}\big), h\big(\{X(s) \cap Q_K(\vp')\}_{s \in \jdm(tT)}\big)\Big) \\
        &\quad= \Cov\Big(h\big(\{X^{0, \vp}(s) \cap Q_K(\vp)\}_{s \in \jdm(0)}\big), h\big(\{X^{tT, \vp'}(s) \cap Q_K(\vp')\}_{s \in \jdm(tT)}\big)\Big).
\end{align*}
Furthermore, by the independence property of the Poisson point process,
$$\Cov\Big(h\big(\{X^{0, \vp}(s) \cap Q_K(\vp)\}_{s \in \jdm(0)}\big), h\big(\{(X^{tT, \vp'} \sm X^{0, \vp})(s) \cap Q_K(\vp')\}_{s \in \jdm(tT)}\big)\Big) = 0.$$
Therefore, it suffices to show that
$\lim_{T\tff}\sup_{\vp, \vp'}\P\big( X^{0, \vp} \cap X^{tT, \vp'} \ne \es \big) = 0.$
Writing $\vp$ and $\vp'$ as finite sets of points, the claim reduces to
$$\lim_{T\tff}\sup_{x, x' \in \R^d}\E\big[\#( X^{0, \{x\}} \cap X^{tT, \{x'\}} ) \big) = 0.$$
But now Lemma \ref{diff_node_lem} gives that 
$\lim_{t \tff} \sup_{x'' \in \R^d}\P(X_0(t) \in Q_K(x'')) = 0$
so that an application of Palm calculus concludes the proof.
\enp

\subsection{Dense regime}
\label{dense_sec}

As announced in Section \ref{out_sec}, the key step to prove Proposition \ref{fr_thm} is the second moment method. To carry out this program, we need decorrelation of the connection intervals at distant time points. Whereas Proposition \ref{decorr_pr} provides such a decorrelation property in a conditional setting, in the dense regime $\a < d/2$ this decorrelation needs to be strengthened into an unconditional result. This will be achieved by applying the law of total covariances, i.e., 
$$\Cov(X, X') = \E[\Cov(X, X' \ba \FF)] + \Cov(\E[X \ba \FF], \E[X' \ba \FF])],$$
for any square-integrable random variables $X, X'$ and $\s$-algebra $\FF$.

%
%
\bep[Proof of Proposition \ref{fr_thm}, $\a < d/2$]
Let $t \le 1$ be arbitrary. We want to show that 
$$\lim_{T\tff}\Cov\big(g\big(\idmm(0)\big), g\big(\idmm(tT)\big) \big)= 0.$$
To achieve this goal, we will apply two times the law of total covariance.

First, by stationarity, similarly as in \eqref{sp_eq}, we can express the conditional expectation given $X_0$ in the form 
$$\E\big[g\big(\idm(tT)\big) \ba X_0\big] = S''\big(\{\one\{X_0(s) - X_0(tT) \}_{s \in \jdm(tT)})\big)$$
for a suitable choice of $S''$. Now, note that the jump times together with the jump directions form an independently marked Poisson point process. Since the intervals $[-M, M]$ and $tT + [-M, M]$ are disjoint for sufficiently large $T$, the independence property of the Poisson point process gives that
$$\Cov\big(\E\big[g\big(\idmm(0)\big)\ba X_0], \E\big[g\big(\idmm(tT)\big)\ba X_0\big] \big) = 0.$$
Hence, by the law of total covariance, it suffices to show that
$$\lim_{T\tff}\E\Big[\Cov\big(g\big(\idmm(0)\big), g\big(\idmm(tT)\big)  \ba X_0\big)\Big] = 0.$$

Combining the law of total covariance with Propositions \ref{decorr_pr} and \ref{xyc_pr} reduces this task to proving that
\begin{align}
	\label{dense_c_eq}
	\lim_{T\tff}\E\Big[\Cov\big(&S'\big(\{X_0(s)\}_{s \in \jdm(0)}, \#Y(0, k)\big),\\
	&	S'\big(\{X_0(s) - X_0(tT)\}_{s \in \jdm(tT)}, \#Y(tT, k)\big)   \ba X_0\big)\Big] = 0.
\end{align}
Now, setting $\a' = (\a + d/2)/2$, in the dense regime the sinks $Y(0, k)$ are contained in $B_{T^{\a'/d}}(o)$ whp, and the sinks $Y(tT, k)$ do not hit $B_{T^{\a'/d}}(o)$ whp. Invoking the independence property of the Poisson process in disjoint domains, this establishes the vanishing of the covariance in \eqref{dense_c_eq}, thereby concluding the proof of Proposition \ref{fr_thm} in the dense regime.
\enp

\subsection{Sparse regime}
\label{sparse_sec}

The high-level proof structure of Proposition \ref{fr_thm} in the sparse regime $\a > d/2$ is similar to that in dense regime discussed in Section \ref{dense_sec}. However, since the sinks $Y$ are sparse in relation to the movement of the typical node they are not affected by the long-time averaging of the movement of the typical node. Therefore, we apply the second moment method conditioned on $Y$.

%
%
\bep[Proof of Proposition \ref{fr_thm}, $\a > d/2$]
As announced above, the proof consists of two steps. First, we replace $\tktb(gh)$ by the conditional expectation $\E[\tktb(gh)\ba Y]$ invoking the conditional second moment method. Second, we derive a more concise representation of the latter expression. The key observation in the sparse regime is that the set of relevant sinks does not change over time. That is, whp, $Y(0, k) = Y(t, k)$ for all $t \le T$.

In order to carry out the second moment method, we note that the jump times together with jump directions form a marked Poisson point process. Since $[-M, M] \cap (tT + [-M, M]) = \es$ for sufficiently large $T > 0$, we deduce that 
\begin{align*}
	\E\Big[\Cov\big(&S\big( \{X_0(s)\}_{s \in \jdm(0)}, \#Y(0, k)\big), \\
	&S\big(\{ X_0(s) - X_0(tT)\}_{s \in \jdm(tT)}, \#Y(tT, k)\big) \ba Y\big)\Big] = 0.
\end{align*}
Hence, combining the law of total covariance with Propositions \ref{decorr_pr} and \ref{xyc_pr} we obtain the asserted decorrelation
$$\lim_{T\tff}\E\Big[\Cov\big(g\big(\idm(0, \Xi^L(Y(0, k)))\big), g\big(\idm(tT, \Xi^L(Y(tT, k)))\big)  \ba Y\big)\Big] = 0.$$

The next step is the identification of the conditional expectation $\E[\tktb(gh) \ba Y]$, i.e., of 
$$\int_0^1\E\big[g\big(\idm(tT, \Xi^L(Y(tT, k)))\big) \ba Y\big]h(t)\d t.$$
Now, by applying Proposition \ref{xyc_pr} and recalling that the relevant sinks do not change over time, the latter integral becomes
\begin{align*}
	&\int_0^1\E\big[S\big(\{X_0(s) - X_0(tT)\}_{s \in \jdm(tT)}, \#Y(0, k)\big) \ba Y\big]h(t)\d t \\
	&\quad=\E\big[S\big(\{X_0(s)\}_{s \in \jdm(0)}, \#Y(0, k)\big) \ba Y\big] \int_0^1h(t)\d t,
\end{align*}
where the identity follows from the time stationarity of the movement model. Moreover, $\#Y(0, k)$ is a Poisson random variable with intensity $\ns$. Thus, inserting the definition of $S$ leads to the limiting representation asserted in Proposition \ref{fr_thm}.
\enp

\subsection{Critical regime}
\label{crit_sec}

Also in the critical regime $\a = d/2$, we follow the blueprint from Sections \ref{dense_sec} and \ref{sparse_sec}. This time, $\s(X_0, Y)$ is the appropriate $\s$-algebra to condition on, so that the conditional decorrelation has already been established in Proposition \ref{decorr_pr}. Conversely, the limit identification is now more involved. 

\bep[Proof of Proposition \ref{fr_thm}, $\a = d/2$]
As explained in the preceding paragraph, due to Propositions \ref{decorr_pr} and \ref{xyc_pr}, the remaining task is the identification of the conditional expectation
$$\int_0^1S\big(\{s \in \jdm(tT)\co X_0(s) - X_0(tT)\}, \#Y(tT, k)\big)h(t)\d t.$$
as $T\tff$.  The first step is to show that we may condition on $X_0(tT)$, i.e., replace the latter expression by 
$$\int_0^1S''(\#Y(tT, k))h(t)\d t,$$
where 
$S''(n) := \E\big[S\big(\{s  \in \jdm(tT)\co X_0(s) \}, n\big)\big].$
To ease notation, we set
$$\bar S_t := S\big(\{X_0(s) - X_0(tT)\}_{s \in \jdm(tT) }, \#Y(tT, k)\big) - S''(\#Y(tT, k)).$$
Since $\E[\bar S_t \ba X_0(tT), Y] = 0$, we deduce that
\begin{align*}
	\E\Big[\Big(\int_0^1\bar S_t h(t)\d t\Big)^2\Big] = \int_0^1\int_0^1\E\big[\Cov\big(\bar S_s, \bar S_t \ba Y\big)\big] h(s) h(t) \d t\d s.
\end{align*}
Again, since $X_0$ is a compound Poisson process, we deduce that 
$\lim_{T\tff}\E\big[\Cov\big(\bar S_s, \bar S_t \ba Y\big)\big] = 0$
for every $s\ne t$, thereby completing the proof of the first step.

Finally, in order to identify the distributional limit of $\int_0^1S''(\#Y(tT, k))h(t)\d t$ as $T \tff$, we proceed along the lines of \cite{mscale}. To render the presentation self-contained, we recall the main steps of the proof. Introducing the unit-intensity process $Y' := Y/\sqrt T$, we can represent $\int_0^1S''(\#Y(tT, k))h(t)\d t$ in the form
$$F(Y', \{X_0(tT)/\sqrt T\}_{t \le 1}) := \int_0^1S''\big(Y'(B_{\nsp}(X_0(tT)/\sqrt T))\big)h(t)\d t.$$
Now, let $f\co [0, 1] \rightarrow [0, 1]$ be an arbitrary Lipschitz function with Lipschitz constant 1. We claim that 

$$\lim_{T\tff}\E\big[f(F(Y', \{X_0(tT)/\sqrt T\}_{t \le 1}))\big] = \E\big[f(F(Y', \{W_t\}_{t \le 1}))\big].$$
By the invariance principle and the continuous mapping theorem, it suffices to show that the map sending a trajectory $\g $ to $\E\big[f(F(Y', \g))\big]$ is continuous outside a zero-set with respect to Brownian motion. Now, let $\{\g_n\}_n$ be a sequence of trajectories. Then, 
\begin{align*}
	\E\big[|f(F(Y', \g_n)) - f(F(Y', \g))|\big] &\le \int_0^1\E\big[\big|S''\big(Y'(B_{\nsp}(\g_n(t)))\big) -S''\big(Y'(B_{\nsp}(\g(t)))\big)\big|\big] \d t\\
	&\le \int_0^1\E\big[Y'\big(B_{\nsp}(\g_n(t))\De B_{\nsp}(\g(t))\big)\big] \d t\\
	&=\int_0^1\big|B_{\nsp}(\g_n(t))\De B_{\nsp}(\g(t))\big|\d t.
\end{align*}
Now, we conclude the proof by noting that the right-hand side tends to 0 as $\g_n \rightarrow \g$ in the sup norm.
\enp

	\section*{Acknowledgements}
The authors thank Quentin Fran{\c c}ois for providing data sets of connection intervals that formed the basis for the simulations presented in the manuscript. This research was supported by Orange S.A., France grant CRE G09292, the German Research Foundation under Germany's Excellence Strategy MATH+: The Berlin Mathematics Research Center, EXC-2046/1 project ID: 390685689, and the Leibniz Association within the Leibniz Junior Research Group on Probabilistic Methods for Dynamic Communication Networks as part of the Leibniz Competition.

\bibliographystyle{abbrv}
\bibliography{refs}
\end{document}